\documentclass{amsart}
\usepackage{graphicx, color}
\usepackage{amscd}
\usepackage{amsmath}
\usepackage{amsfonts}
\usepackage{amssymb}
\usepackage{mathrsfs} 
\newtheorem{theorem}{Theorem}
\newtheorem{lemma}[theorem]{Lemma}

\newtheorem{proposition}[theorem]{Proposition}

\newtheorem{corollary}[theorem]{Corollary}

\numberwithin{equation}{section}
\newcommand{\erre}{\mathbb R}
\newcommand{\R}{\mathbb R}
\newcommand{\N}{\mathbb N}
\newcommand{\RR}{\mathbb R}
\newcommand{\NN}{\mathbb N}

\begin{document}

\title[Some hemivariational inequalities]{Some hemivariational inequalities in the Euclidean space}

\author{Giovanni Molica Bisci}
\address[G. Molica Bisci]{Dipartimento di Scienze Pure e Applicate (DiSPeA) - 
Universit\`{a} degli Studi di Urbino Carlo Bo, 
Piazza della Repubblica 13,
61029 Urbino, Italy}
\email{\tt giovanni.molicabisci@uniurb.it}

\author{Du\v{s}an D. Repov\v{s}}
\address[Du\v{s}an D. Repov\v{s}]{Faculty of Education, 
University of Ljubljana, Kardeljeva pl. 16, 
1000 Ljubljana, Slovenia
\&
Faculty of Mathematics and Physics,
University of Ljubljana,
Jadranska 21,
1000 Ljubljana, Slovenia \& 
Institute of Mathematics, Physics and Mechanics, 
Jadranska 19,
1000 Ljubljana, Slovenia}
\email{\tt dusan.repovs@guest.arnes.si}


\begin{abstract}
The purpose of this paper is to study the existence of weak solutions for some classes of
hemivariational problems in the Euclidean space $\R^d$ ($d\geq 3$). These hemivariational inequalities have a
variational structure and, thanks to this, we are able to find a non-trivial weak solution for them
by using variational methods and a non-smooth version of the Palais principle of symmetric criticality for locally Lipschitz continuous functionals, due to Krawcewicz and Marzantowicz.
   The main tools in our approach are based on appropriate theoretical arguments on suitable subgroups of the orthogonal group $O(d)$ and their actions on the Sobolev space $H^1(\R^d)$. Moreover, under an additional hypotheses on the dimension $d$ and in the presence of symmetry on the nonlinear datum, the existence of multiple pairs of sign-changing solutions with different symmetries structure has been proved. In connection to classical Schr\"{o}dinger equations a concrete and meaningful example of an application is presented.
\end{abstract}

\maketitle

\tableofcontents

\section{Introduction}\label{sec:introduzione}

The aim of this paper is to study some nonlinear eigenvalue problems for certain classes of hemivariational inequalities
that depend on a real parameter. For instance, the
motivation for such a study comes from the investigation of
perturbations, usually determined in terms of parameters. The hemivariational
inequalities appears as a generalization of the variational inequalities and their study is based on the notion of Clarke subdifferential of a
locally Lipschitz function. The theory of hemivariational inequalities appears as a new
field of Non-smooth Analysis; see \cite[Part I - Chapter II]{KRO1} and the references therein.

\indent More precisely, we study the following hemivariational inequality problem:\par
\smallskip

\begin{itemize}
\item[$({S}_\lambda)$] {\it Find $u\in H^1(\R^d)$ such that
\begin{equation*}
\left\{\begin{array}{l} \displaystyle\int_{\R^d}\nabla u(x)\cdot \nabla \varphi(x)dx +\int_{\R^d}u(x)\varphi(x))dx\\
\qquad \qquad\qquad \qquad\quad\quad +\displaystyle\lambda\int_{\R^d}W(x)F^0(u(x);-\varphi(x))dx\geq 0,\,\,\,\,\,\,\\
\forall\,\varphi\in H^1(\R^d).
\end{array}\right.
\end{equation*}}
\end{itemize}
\smallskip
 Here $(\R^d,|\cdot|)$ denotes the Euclidean space (with $d\geq 3$), $F:\R \rightarrow \R$ is a locally Lipschitz continuous function, whereas
 $$
 F^0(s;z):=\limsup_{\substack{y\to s\\ t\to 0^+}}{{F(y+tz)-F(y)}\over t}
 $$
 is the generalized directional derivative of $F$ at the point $s\in \R$ in
the direction $z\in \R$; see the classical monograph of Clarke \cite{CLARKE2} for details. Finally, $W\in L^{\infty}(\R^d)\cap L^{1}(\R^d)\setminus\{0\}$ is a non-negative radially symmetric map and $\lambda$ is a positive real parameter.

We assume that there exist $\kappa_1>0$ and $q\in (2,2^*)$, where $2^*=2d/(d-2)$, such that
\begin{equation}\label{crescita}
|\zeta|\leq \kappa_1 (1+|s|^{q-1}),\quad \forall \zeta\in \partial F(s),\quad \mbox{for every}\, s\in \erre,
\end{equation}
where $\partial F(s)$ denotes the generalized gradient of the function $F$ at $s\in \R$ (see Section \ref{sec:preparatory}).\par
With the above notations the main result reads as follows.
\begin{theorem}\label{Main1}
 Let $F:{\R}\to{\R}$ be a locally Lipschitz continuous function with $F(0)=0$ and satisfying the growth condition \eqref{crescita} for some $q\in (2,2^*)$, in addition to
\begin{equation}\label{azero}
\limsup_{s\rightarrow0^+} \frac{F(s)}{s^2}=+\infty
\,\,\,\,\,\,\,{and}\,\,\,\,\,\,\liminf_{s\rightarrow0^+} \frac{F(s)}{s^2}>-\infty.
\end{equation}
Moreover, let $W\in L^{\infty}(\R^d)\cap L^{1}(\R^d)\setminus\{0\}$ be a non-negative radially symmetric map. Then the following facts hold:

\begin{itemize}
\item[$(a_1)$] There exists a positive number $\lambda^{\star}$
such that, for every $\lambda\in (0,\lambda^{\star})$,
the problem $(S_\lambda)$ admits at least one non-trivial radial weak solution $u_{\lambda}\in H^1(\R^d)$ with
$|u_\lambda(x)|\rightarrow 0$ as $|x|\rightarrow \infty$.

\item[$(a_2)$] If $d> 3$ and $F$ is even then there exists a positive number $\lambda_{\star}$
such that for every $\lambda\in (0,\lambda_{\star})$,
the problem $(S_\lambda)$ admits at least
 $$
 \zeta^{(d)}_S:=1+(-1)^{d}+\left[\frac{d-3}{2}\right]
 $$
 pairs of non-trivial weak solutions $\{\pm u_{\lambda,i}\}_{i\in J'_d}\subset H^1(\R^d)$ with $|u_{\lambda,i}(x)|\rightarrow 0$, as $|x|\rightarrow \infty$, for every $i\in J'_d:=\{1,...,\zeta^{(d)}_S\}$, and with different symmetries structure. More precisely, if $d\neq 5$ problem $(S_\lambda)$ admits at least
 $$\tau_d:= \zeta^{(d)}_S-1$$ pairs of sign-changing weak solutions.
\end{itemize}
\end{theorem}

Here, the symbol $[\cdot]$ denotes
the integer function.\par

\smallskip

The proof of the above result is based on variational method in the nonsmooth setting. As it is well known, the lack of a compact embeddings of the Sobolev space $H^1(\R^d)$ into Lebesgue spaces produces several difficulties for exploiting variational methods. In order to recover compactness, the first task is to construct certain subspaces of $H^1(\R^d)$ containing invariant
functions under special actions defined by means of carefully chosen subgroups of
the orthogonal group $O(d)$. Subsequently, a locally Lipschitz continuous function
is constructed which is invariant under the action of suitable subgroups of $O(d)$,
whose restriction to the appropriate subspace of invariant functions
admits critical points.\par

 Thanks to a nonsmooth version of the principle of symmetric criticality
obtained by Krawcewicz and Marzantowicz \cite{KaMa}, these points will also be critical points of the original
functional, and they are exactly weak solutions
of problem $(S_\lambda)$. The abstract critical point result that we employ here is a nonsmooth version of the variational principle established by Ricceri \cite{ricceri}; see Bonanno and Molica Bisci \cite{BoMo} for details.

Moreover, we also emphasize that the multiplicity property stated in Theorem \ref{Main1} - part $(a_2)$ is obtained by using the group-theoretical approach developed by Krist\'{a}ly, Moro\c{s}anu, and O'Regan \cite{KRO}; see Subsection \ref{Zzz}. Thanks to this analysis, we are able to construct
$$
 \zeta^{(d)}_S:=1+(-1)^{d}+\left[\frac{d-3}{2}\right]
 $$
subspaces of $H^1(\R^d)$ with different symmetries properties. In addition, when $d\neq 5$, there are
$$
\tau_d:=(-1)^{d}+\left[\frac{d-3}{2}\right]
$$
of these subspaces which do not contain radial symmetric functions; see the quoted paper \cite{bw} due to Bartsch and Willem, as well as \cite[Theorem 2.2]{KRO}.

We point out that some almost straightforward computations in \cite{Molica} are
adapted here to the non-smooth case. However, due to the non-smooth framework, our abstract procedure, as well as the setting of the main results,
is different from the results contained in \cite{Molica}, where the continuous case was studied; see Section \ref{sec finale} for additional comments and remarks.

\indent The manuscript is organized as follows. In Section~\ref{sec:preparatory}
we set some notations and recall some properties of the functional space we shall work in.
In order to apply critical point methods to problem $(S_\lambda)$, we need to work in a subspace of the functional space $H^1(\R^d)$
 in particular, we give some tools which will be useful in the paper (see Propositions \ref{derivata0} and Lemma \ref{semicontinuity}). In Section~\ref{sec:Main} we study problem~$(S_\lambda)$ and we prove our existence result (see Theorem \ref{Main1}). Finally, we study the existence of multiple non-radial solutions to the problem~$(S_\lambda)$ for $\lambda$ sufficiently small. In connection to classical Schr\"{o}dinger equations in the continuous setting (see, among others, the papers \cite{BLW,BPW,bw2, BLions}) a meaningful example of an application is given in the last section.\par
 We refer to the books \cite{AM2, KRO1,Wi} as general references on the subject treated in the paper.

\section{Abstract framework}\label{sec:preparatory}

\indent Let $(X,\Vert\cdot\Vert_X)$ be a real Banach space.
We denote by $X^*$ the dual space of $X$, whereas
$\langle\cdot,\cdot\rangle$ denotes the duality pairing between
$X^*$ and $X$.\par
 A function $J:X\to\erre$ is called \textit{locally
Lipschitz continuous} if to every $y\in X$ there correspond a
neighborhood $V_y$ of $y$ and a constant $L_y\geq 0$ such that
$$\vert J(z)-J(w)\vert\leq L_y\Vert z-w\Vert_X,\quad(\forall\, z,w\in V_y).$$
\indent If $y,z\in X$, we write $J^{0}(y;z)$ for the generalized
directional derivative of $J$ at the point $y$ along the
direction $z$, i.e.,
$$J^{0}(y;z):=\limsup_{\substack{w\to y\\ t\to 0^+}}{{J(w+tz)-J(w)}\over t}.$$
\indent The generalized gradient of the function $J$ at $y\in X$,
denoted by $\partial J(y)$, is the set
$$\partial J(y):=\left\{ y^*\in X^*:\, \langle y^*,z\rangle\leq J^{0}(y;z),
\;\forall\, z\in X\right\}.$$
 \indent The basic properties of generalized directional derivative and generalized gradient which we shall use here were studied in \cite{CHANG,CLARKE2}.\par

 The following lemma displays some useful properties of the notions
introduced above.

\begin{lemma}\label{gd}
If $I,J:X\to\mathbb{R}$ are locally Lipschitz continuous functionals, then
\begin{itemize}
\item[$(i)$]
$J^0(y;\cdot)$ is positively homogeneous,
sub-additive, and continuous for every $y\in X;$

\item[$(ii)$] $J^0(y;z)=\displaystyle\max\{\langle y^*,z\rangle:{y^*\in \partial J(z)}\}$ for every $y,z\in X;$

\item[$(iii)$]
$J^0(y;-z)=(-J)^0(y;z)$ for every $y,z\in X;$

\item[$(iv)$]
 if $J\in C^1(X)$, then $J^0(y;z)=\langle J'(y),z\rangle$
for every $y,z\in X;$

\item[$(v)$] 
$(I+J)^0(y;z)\leq I^0(y;z)+J^0(y;z)$ for every
 $y,z\in X$. Moreover, if $J$ is is continuously
G\^{a}teaux differentiable, then $(I+J)^0(y;z)=I^0(y;z)+J'(y;z)$ for every
 $y,z\in X$.
\end{itemize}
\end{lemma}

See \cite{GP} for details.\par

\smallskip

Further, a point $y\in X$ is called a (generalized)
\textit{critical point} of the locally Lipschitz continuous function $J$ if $0_{X^*}\in \partial J(y)$,
 i.e.
$$J^0(y;z)\geq 0,$$
 for every $z\in X$.\par
  Clearly, if $J$ is a continuously G\^{a}teaux differentiable at $y\in X$, then $y$
becomes a (classical) critical point of $J$, that is $J '(y) = 0_{X^*}$.

\indent For an exhaustive overview of the non-smooth calculus we refer to the
monographs \cite{CHANG, CLARKE2, Mot1, Mot2}. Further, we cite the book \cite{KRO1} as a general reference on this subject.

To make the nonlinear methods work, some
careful analysis of the fractional spaces involved is necessary.
Assume $d\geq 3$ and let $H^1(\R^d)$ be the standard Sobolev space endowed by the inner product
$$
\langle u, v\rangle:=\int_{\R^d}\nabla u(x)\cdot \nabla v(x)dx +\int_{\R^d}u(x)v(x)dx,\quad\forall\, u,v\in H^1(\R^d)
$$
and the induced norm
$$
\|u\|:=\left(\int_{\R^d}|\nabla u(x)|^2dx+\int_{\R^d}|u(x)|^2dx\right)^{1/2},
$$
for every $u\in H^1(\R^d)$.\par

\indent In order to prove Theorem \ref{Main1} we apply the principle of symmetric criticality together with the following critical point theorem proved in \cite{BoMo} by Bonanno and Molica Bisci.

\begin{theorem}\label{BMB}
Let $X$ be a reflexive real Banach space and let $\Phi,\Psi:X\to\R$
be locally Lipschitz continuous functionals such that $\Phi$ is
sequentially weakly lower semicontinuous and coercive. Furthermore, assume that $\Psi$
is sequentially weakly upper semicontinuous. For every $r>\inf_X
\Phi$, put
$$
\varphi(r):=\inf_{u\in\Phi^{-1}((-\infty,r))}\frac{\displaystyle\left(\sup_{v\in\Phi^{-1}((-\infty,r))}\Psi(v)\right)-\Psi(u)}{r-\Phi(u)}.
$$
Then for each $r>\inf_X\Phi$ and each
$\lambda\in\left]0,{1}/{\varphi(r)}\right[$, the restriction
of $\mathcal J_\lambda:=\Phi-\lambda\Psi$ to
$\Phi^{-1}((-\infty,r))$ admits a global minimum, which is a
critical point $($local minimum$)$ of $\mathcal J_\lambda$ in $X$.
\end{theorem}

The above result represents a nonsmooth version of a variational principle established by Ricceri in \cite{ricceri}.
\smallskip

For completeness, we also recall here the principle of
symmetric criticality of Krawcewicz and Marzantowicz which represents a non-smooth version of the celebrated
result proved by Palais in \cite{pala}. We point out that the result proved in \cite{KaMa}
was established for sufficiently smooth Banach $G$-manifolds.
We will use here a particular form of this result that is valid for Banach spaces.

 An \textit{action} of a compact Lie group $G$ on the Banach space $(X,\|\cdot\|_X)$ is a continuous map
$$
*: G\times X\rightarrow X: (g,y)\mapsto g*y,
$$
\noindent such that
$$
1*y=y,\,\,\,
(gh)*y=g*(h*y),\,\,\, y\mapsto g*y\,\,\,{\rm is\,linear}.
$$
\indent The action $*$ is said to be \textit{isometric} if $\|g*y\|_X=\|y\|_X$, for every $g\in G$ and $y\in X$. Moreover,
 the space of $G$-invariant points is defined by
$$
Fix_G(X):=\{y\in X:g*y=y,\forall g\in G\},
$$
and a map $h:X\rightarrow \R$ is said to be $G$-\textit{invariant} on $X$ if
$$
h(g*y)=h(y)
$$
for every $g\in G$ and $y\in X$.\par

\begin{theorem}\label{Palais}
Let $X$ be a Banach space, let $G$ be a compact topological group acting linearly and
isometrically on $X$, and $J: X\rightarrow\R$ a locally Lipschitz, $G$-invariant functional. Then every
critical point of $\mathcal{J}: Fix_G(X)\rightarrow\R$ is also a critical point of $J$.
\end{theorem}

For details see, for instance, the book \cite[Part I - Chapter 1]{KRO1} and Krawcewicz and Marzantowicz \cite{KaMa}.\par

\subsection{Group-theoretical arguments}\label{Zzz}

Let $O(d)$ be the orthogonal group in $\R^d$ and let $G\subseteq O(d)$ be a subgroup. Assume that $G$ acts on the space $H^1(\R^d)$. Hence, the set of fixed points of $H^1(\R^d)$, with respect to $G$, is clearly given by
$$
Fix_{G}(H^1(\R^d)):=\{u\in H^1(\R^d):gu=u,\forall g\in G\}.
$$
\noindent We note that, if $G=O(d)$ and the action is the standard linear isometric map defined by
$$
gu(x):=u(g^{-1}x),\quad\forall\, x\in \R^d\quad {\rm and}\quad g\in O(d)
$$
then $Fix_{O(d)}(H^1(\R^d))$ is exactly the \textit{subspace of radially symmetric functions} of $H^1(\R^d)$, also denoted by $H^1_{\rm rad}(\R^d)$. Moreover, the following embedding
\begin{equation}\label{compattezza}
Fix_{O(d)}(H^1(\R^d))\hookrightarrow L^{q}(\R^d)
\end{equation}
is continuous (resp. compact), for every $q\in [2,2^*]$ (resp. $q\in (2,2^*)$). See, for instance, the celebrated paper \cite{Lions}.\par

Let either $d=4$ or $d\geq 6$ and consider the subgroup $H_{d,i}\subset O(d)$ given by
\[
H_{d,i}:= \left\{
\begin{array}{ll}
\displaystyle O(d/2)\times O(d/2) & \mbox{ if $i=\displaystyle\frac{d-2}{2}$} \\
\displaystyle O(i+1)\times O(d-2i-2)\times O(i+1) & \mbox{ if $i\neq\displaystyle\frac{d-2}{2}$},
\end{array}
\right.
\]
for every $i\in J_d:=\{1,...,\tau_d\}$, where
$$
\tau_d:=(-1)^{d}+\displaystyle\left[\frac{d-3}{2}\right].
$$
\indent Let us define the involution $\eta_{_{H_{d,i}}}:\R^d\rightarrow \R^d$ as follows
\[
\eta_{_{H_{d,i}}}(x):= \left\{
\begin{array}{ll}
\displaystyle (x_3,x_1) & \mbox{ if $i=\displaystyle\frac{d-2}{2}$ and $x:=(x_1,x_3)\in \R^{d/2}\times\R^{d/2}$} \\
\displaystyle (x_3,x_2,x_1) & \mbox{ if $i\neq\displaystyle\frac{d-2}{2}$ and $x:=(x_1,x_2,x_3)\in \R^{i+1}\times \R^{d-2i-2}\times\R^{i+1}$},
\end{array}
\right.
\]
for every $i\in J_d$.\par
By definition, one has $\eta_{_{H_{d,i}}}\notin H_{d,i}$, as well as
$$
\eta_{_{H_{d,i}}} H_{d,i}\eta_{_{H_{d,i}}}^{-1}=H_{d,i},\quad {\rm and}\quad\eta_{_{H_{d,i}}}^2=\textrm{id}_{\R^d},
$$
for every $i\in J_d$.\par
\indent Moreover, for every $i\in J_d$, let us consider the compact group
$$
H_{d,\eta_i}:=\langle H_{d,i}, \eta_{_{H_{d,i}}}\rangle,
$$
that is $H_{d,\eta_i}=H_{d,i}\cup \eta_{_{H_{d,i}}} H_{d,i}$, and the action $\circledast_i:H_{d,\eta_i}\times H^1(\R^d)\rightarrow H^1(\R^d)$ of $H_{d,\eta_i}$ on $H^1(\R^d)$ given by
\begin{equation}\label{isometr}
h \circledast_i u(x):= \left\{
\begin{array}{ll}
\displaystyle u(h^{-1}x) & \mbox{ if $h\in H_{d,i}$} \\
\displaystyle -u(g^{-1}\eta_{_{H_{d,i}}}^{-1}x) & \mbox{ if $h=\eta_{_{H_{d,i}}} g\in H_{d,\eta_i}\setminus H_{d,i}$},\, g\in H_{d,i}
\end{array}
\right.
\end{equation}
for every $x\in \R^d$.\par
We note that $\circledast_i$ is defined for every element of $H_{d,\eta_i}$. Indeed, if $h\in H_{d,\eta_i}$, then either $h\in H_{d,i}$ or $h=\tau g\in H_{d,\eta_i}\setminus H_{d,i}$, with $g\in H_{d,i}$.
Moreover, set
$$
Fix_{H_{d,\eta_i}}(H^1(\R^d)):=\{u\in H^1(\R^d):h\circledast_i u=u,\forall h\in H_{d,\eta_i}\},
$$
for every $i\in J_d$.\par
\indent Following Bartsch and Willem \cite{bw}, for every $i\in J_d$, the embedding
\begin{equation}\label{compattezza2}
Fix_{H_{d,\eta_i}}(H^1(\R^d))\hookrightarrow L^{q}(\R^d)
\end{equation}
is compact, for every $q\in (2,2^*)$.

 \begin{proposition}\label{geometry}
 With the above notations, the following properties hold:
 \begin{itemize}
 \item[] if $d=4$ or $d\geq 6$, then
\begin{equation}\label{intersezione}
Fix_{H_{d,\eta_i}}(H^1(\R^d))\cap Fix_{O(d)}(H^1(\R^d))=\{0\},
\end{equation}
for every $i\in J_d$;
\end{itemize}
\begin{itemize}
 \item[] if $d=6$ or $d\geq 8$, then
\begin{equation}\label{intersezione2}
Fix_{H_{d,\eta_i}}(H^1(\R^d))\cap Fix_{H_{d,\eta_j}}(H^1(\R^d))=\{0\},
\end{equation}
for every $i,j\in J_d$ and $i\neq j$.
\end{itemize}
\end{proposition}
See \cite[Theorem 2.2]{KRO} for details.\par

\noindent From now on, for every $u\in L^\ell(\R^d)$ and $\ell\in [2,2^*)$, we shall denote
$$\|u\|_\ell:=\left(\displaystyle\int_{\R^d}|u(x)|^\ell dx\right)^{1/\ell},$$  and
$$
\|W\|_\infty:=\displaystyle{\rm esssup}_{x\in \R^d} |W(x)|,\quad\quad\|u\|_p:=\left(\displaystyle\int_{\R^d}|u(x)|^p dx\right)^{1/p},
$$
for every $p\in [2,2^*)$.\par
Moreover, let $\Psi:H^1(\R^d)\rightarrow \erre$ given by
$$
\Psi(u):=\int_{\R^d}W(x)F(u(x))dx,\quad \forall\, u\in H^1(\R^d).
$$
The following locally Lipschitz property holds.

\begin{lemma}\label{loclip}
Assume that
condition \eqref{crescita} holds for some $q\in \left(2,2^*\right)$ and $F(0)=0$. Furthermore, let $W\in L^{\infty}(\R^d)\cap L^{1}(\R^d)\setminus\{0\}$. Then the extended functional $\Psi^{e}:L^{q}(\R^d)\rightarrow \erre$ defined by
$$
\Psi^{e}(u):=\int_{\R^d}W(x)F(u(x))dx,\quad \forall\, u\in L^{q}(\R^d)
$$
is well-defined and locally Lipschitz continuous on $L^{q}(\R^d)$.
\end{lemma}
{\it Proof.}
It is clear that $\Psi^{e}$ is well-defined. Indeed, by using Lebourg's mean value theorem, fixing $t_1,t_2\in \R$, there exist
$\theta\in (0,1)$ and $\zeta_\theta\in \partial F(\theta t_1+(1-\theta)t_2)$ such that
\begin{eqnarray}\label{Lebou}
F(t_1)-F(t_2)=\zeta_\theta (t_1-t_2).
\end{eqnarray}
Since $F(0)=0$, by using \eqref{Lebou} and condition \eqref{crescita}, our assumptions on $W$ and the H\"{o}lder inequality gives that
\begin{eqnarray}\label{ineq}
  \int_{\mathbb R^d} W(x)F(u(x))dx &\leq & \kappa_1\left(\int_{\mathbb R^d}W(x)|u(x)|dx+\int_{\mathbb R^d}W(x)|u(x)|^{q}dx\right)\nonumber\\
   &\leq& \kappa_1\left(\int_{\mathbb R^d}|W(x)|^{\frac{q}{q-1}}dx\right)^{\frac{q-1}{q}}\left(\int_{\mathbb R^d}|u(x)|^qdx\right)^{1/q}\\
   &&+\kappa_1\|W\|_{\infty}\int_{\mathbb R^d}|u(x)|^qdx,\nonumber
\end{eqnarray}
for every $u\in L^{q}(\R^d)$. Hence, inequality \eqref{ineq} yields
\begin{equation}\label{ineq2}
\Psi^{e}(u) \leq \kappa_1\left(\|W\|_{\frac{q}{q-1}}\|u\|_q+\|W\|_{\infty}\|u\|_q^q\right)<+\infty,
\end{equation}
for every $u\in L^{q}(\R^d)$.\par
In order to prove that $\Psi^e$ is locally Lipschitz continuous on $L^{q}(\R^d)$ it is
straightforward to establish that the functional $\Psi^e$ is in fact Lipschitz continuous on $L^{q}(\R^d)$. Now, for a fixed number $r>0$ and arbitrary elements $u,v\in L^{q}(\R^d)$ with $\max\{\|u\|_q,\|v\|_q\}\leq r$, the following estimate holds
\begin{eqnarray}\label{norms24444}
  |\Psi^e(u)-\Psi^e(v)| &\leq& \int_{\R^d}W(x)\left|F(u(x))-F(v(x))\right|dx\nonumber\\
   &\leq& \kappa_1\int_{\mathbb R^d}W(x)\left(1+|u(x)|^{q-1}+|v(x)|^{q-1}\right)|u(x)-v(x)|dx \\
    &\leq& \kappa_1(\|W\|_{\frac{q}{q-1}}\|u-v\|_q+\|W\|_{\infty}(\|u\|_{q}^{q-1}+\|v\|_{q}^{q-1})\|u-v\|_q)\nonumber\\
   &\leq& \kappa_2\|u-v\|_q,\nonumber
\end{eqnarray}
where the Lipschitz constant $\kappa_2:=2^{q-2}(\|W\|_{\frac{q}{q-1}}+2r^{q-1}\|W\|_{\infty})\kappa_1$ depends on $r$.

The above inequalities have
been derived by using \eqref{Lebou}, assumption \eqref{crescita} and H\"{o}lder's inequality.
The Lipschitz property on bounded sets for $\Psi^e$ is thus verified.  \hfill$\Box$
\smallskip

\indent A meaningful consequence of the above lemma is the following semicontinuity property.

\begin{corollary}\label{semicontinuity}

Assume that
condition \eqref{crescita} holds for some $q\in \left(2,2^*\right)$ and let $W\in L^{\infty}(\R^d)\cap L^{1}(\R^d)\setminus\{0\}$. Then for every $\lambda>0,$ the functional
 $$
 u\mapsto \frac{1}{2}\|u\|^2-\lambda\Psi|_{Fix_{Y}(H^1(\R^d))}(u),\quad\forall\, u\in Fix_{Y}(H^1(\R^d))
 $$
 is sequentially weakly lower semicontinuous on $Fix_{Y}(H^1(\R^d))$, where either $Y=O(d)$ or $Y=H_{d,\eta_i}$ for some $i\in J_d$.
\end{corollary}
{\it Proof.}  First, on account of Br\'ezis \cite[Corollaire
III.8]{brezis}, the functional
$u\mapsto\|u\|^2/2$ is sequentially weakly lower semicontinuous on $Fix_{Y}(H^1(\R^d))$. Now, we prove that $\Psi|_{Fix_{Y}(H^1(\R^d))}$ is sequentially weakly continuous. Indeed, let $\{u_j\}_{j\in \mathbb N}\subset Fix_{Y}(H^1(\R^d))$ be a sequence which weakly converges to an element $u_0\in Fix_{Y}(H^1(\R^d))$. Since $Y$ is compactly embedded in $L^q(\R^d)$, for every $q\in (2,2^*)$, passing to a subsequence if necessary, one has $\|u_j-u_0\|_q\rightarrow 0$ as $j\rightarrow \infty$. According to Lemma \ref{loclip}, the extension of $\Psi$ to $L^q(\R^d)$ is locally Lipschitz continuous. Hence, there exists a constant $L_{u_0}\geq 0$ such that
\begin{equation}\label{LocLip}
|\Psi(u_j)-\Psi(u_0)|\leq L_{u_0}\|u_j-u_0\|_q,
\end{equation}
for every $j\in \N$. Passing to the limit in \eqref{LocLip}, we conclude that $\Psi$ is sequentially weakly continuous on $Fix_{Y}(H^1(\R^d))$. The proof is now complete.
\hfill$\Box$

\smallskip

The next result will be crucial in the sequel; see \cite{CLARKE2, Kri, Kri2, Mot1} for related results.

\begin{proposition}\label{derivata0}
Assume that
condition \eqref{crescita} holds for some $q\in \left(2,2^*\right)$ and let $W\in L^{\infty}(\R^d)\cap L^{1}(\R^d)\setminus\{0\}$.
Furthermore, let $E$ be a closed subspace of $H^1(\R^d)$ and denote by $\Psi_E$ the restriction of $\Psi$ to $E$.
Then the following inequality holds
\begin{equation}\label{Derivate}
\Psi^0_E(u;v)\leq \int_{\R^d}W(x)F^0(u(x);v(x))dx,
\end{equation}
for every $u,v\in E$.
\end{proposition}
{\it Proof.} The map $x\mapsto W(x)F^0(u(x);v(x))$ is measurable on $\R^d$. Indeed, $W\in L^{\infty}(\R^d)$ and the function $x\mapsto F^0(u(x);v(x))$  is measurable as the countable
limsup of measurable functions, see p. 16 of \cite{Mot1} for details. Moreover, condition \eqref{crescita} ensures that
$$
\int_{\R^d}W(x)F^0(u(x);v(x))dx<\infty.
$$
Thus the map $x\mapsto W(x)F^0(u(x);v(x))$ belongs to $L^1(\R^d)$.

 The next task is to prove \eqref{Derivate}. To this goal, since $E$ is separable, let us notice that there exist two sequences $\{t_j\}_{j\in \N}\in \R$
and $\{w_j\}_{j\in \N}\subset E$ such that $t_j\rightarrow 0^+$, $\|w_j-u\|\rightarrow 0$ in $E$ and
$$
 \Psi^0_E(u;v)=\displaystyle\lim_{j\rightarrow \infty}{{\displaystyle\Psi_E(w_j+t_jv)-\displaystyle\Psi_E(w_j)}\over t_j}.
 $$
Without loss of generality we can also suppose that $w_j(x)\rightarrow u(x)$ a.e. in $\R^d$ as $j\rightarrow \infty$.

Now, for every $j\in \mathbb{N}$, let us consider the measurable and non-negative function $g_j:\R^d\rightarrow\R\cup \{+\infty\}$ defined by
$$
g_j(x):=\kappa_1|v(x)|(1+|w_j(x)+t_jv(x)|^{q-1}+|w_j(x)|^{q-1})
$$
$$
\quad-{{\displaystyle F(w_j(x)+t_jv(x))-\displaystyle F(w_j(x))}\over t_j},
$$
for a.e. $x\in \R^d$.
Set
$$
I:=\limsup_{j\rightarrow \infty}\left(-\int_{\R^d}W(x)g_j(x)dx\right).
$$
The inverse Fatou's Lemma applied to the sequences $\{Wg_j\}_{j\in \mathbb{N}}$ yields
\begin{equation}\label{1}
I\leq J:=\int_{\R^d} W(x)\limsup_{j\rightarrow \infty}(\alpha_j(x)-\beta_j(x))dx,
\end{equation}
where
$$
\alpha_j(x)={{\displaystyle F(w_j(x)+t_jv(x))-\displaystyle F(w_j(x))}\over t_j},
$$
and
$$
\beta_j(x):=\kappa_1|v(x)|(1+|w_j(x)+t_jv(x)|^{q-1}+|w_j(x)|^{q-1})
$$
for every $j\in \mathbb{N}$ and a.e. $x\in \R^d$.

By setting
$$
\gamma_j:=\int_{\R^d}W(x)\beta_j(x)dx,
$$
one has
\begin{equation}\label{Decompo}
I=\limsup_{j\rightarrow \infty}\left(\int_{\R^d}W(x)\alpha_j(x)dx-\gamma_j\right).
\end{equation}

Now, it is easily seen that there exists a function $k\in L^1(\R^d)$ such that
$$
|\beta_j(x)|\leq k(x),
$$
and $$\beta_j(x)\rightarrow \kappa_1|v(x)|(1+2|u(x)|^{q-1})$$
for a.e. $x\in \R^d$.

Consequently, the Lebesgue's Dominated Convergence Theorem implies that
\begin{equation}\label{limite}
\lim_{j\rightarrow \infty}\gamma_j=\kappa_1\int_{\R^d}W(x)|v(x)|(1+2|u(x)|^{q-1})dx.
\end{equation}
By \eqref{Decompo} and \eqref{limite} it follows that
\begin{eqnarray}\label{Decompo2}
& I=\displaystyle\limsup_{j\rightarrow \infty}{{\displaystyle\Psi_E(w_j+t_jv)-\displaystyle\Psi_E(w_j)}\over t_j}-\lim_{j\rightarrow \infty}\gamma_j\\
&\qquad\qquad=\Psi^0_E(u;v)-\displaystyle\kappa_1\int_{\R^d}W(x)|v(x)|(1+2|u(x)|^{q-1})dx.\nonumber
\end{eqnarray}
\indent Now
\begin{eqnarray}\label{Decompo3}
J\leq J_\alpha-\kappa_1\int_{\R^d}W(x)|v(x)|(1+2|u(x)|^{q-1})dx.
\end{eqnarray}
where
$$
J_\alpha:=\int_{\R^d} W(x)\limsup_{j\rightarrow \infty}\alpha_j(x)dx.
$$
Inequality \eqref{1} in addition to \eqref{Decompo2} and \eqref{Decompo3} yield
\begin{eqnarray}\label{Decompo4}
\Psi^0_E(u;v)\leq J_\alpha.
\end{eqnarray}

Finally,
\begin{eqnarray}\label{FinaleXX}
& J_\alpha=\displaystyle\int_{\R^d} W(x)\limsup_{j\rightarrow \infty}{{\displaystyle F(w_j(x)+t_jv(x))-\displaystyle F(w_j(x))}\over t_j}dx\nonumber\\
&\leq \displaystyle\int_{\R^d}W(x)\lim_{j\rightarrow \infty}{{\displaystyle F(w_j+t_jv)-\displaystyle F(w_j)}\over t_j}dx\\
& \leq \displaystyle \int_{\R^d}W(x)F^0(u(x);v(x))dx.\nonumber
\end{eqnarray}

By \eqref{Decompo4} and \eqref{FinaleXX}, inequality \eqref{Derivate} now immediately follows. \hfill$\Box$

\smallskip

The next result is a direct and easy consequence of Proposition \ref{derivata0}.

\begin{proposition}\label{derivata}
Assume that
condition \eqref{crescita} holds for some $q\in \left(2,2^*\right)$ and let $W\in L^{\infty}(\R^d)\cap L^{1}(\R^d)\setminus\{0\}$.
Let $J_\lambda:H^1(\R^d)\rightarrow \R$ be the functional defined by
$$
J_\lambda(u):=\frac{1}{2}\|u\|^2-\lambda\Psi(u),\quad  \forall\,u\in H^1(\R^d).
$$
Then the functional is locally Lipschitz continuous and its critical points solve $(S_\lambda)$.
\end{proposition}
{\it Proof.} The functional $J_\lambda$ is locally Lipschitz continuous. Indeed, $J_\lambda$ is the sum of the $C^1(H^1(\R^d))$ functional $u\mapsto \|u\|^2/2$ and of the locally Lipschitz continuous functional $\Psi$, see Lemma \ref{loclip}. Now, every critical point of $J_\lambda$ is a weak solution of problem $(S_\lambda)$. Indeed, if $u_0\in H^1(\R^d)$ is a critical point of $J_\lambda$, a direct application of inequality \eqref{Derivate} in Proposition \ref{derivata0} yields
\begin{eqnarray}\label{critical}
& 0\leq J^0_\lambda(u_0;\varphi) = \langle u_0,\varphi\rangle+\lambda(-\Psi)^0(u_0;\varphi) \nonumber \\
&\qquad\qquad\qquad= \langle u_0,\varphi\rangle+\lambda(-\Psi)^0(u_0;\varphi)  \\
&\quad\qquad\leq \langle u_0,\varphi\rangle+\lambda \displaystyle\int_{\R^d}W(x)F^0(u_0(x);-\varphi(x))dx,\nonumber
\end{eqnarray}
for every $\varphi\in H^1(\R^d)$. Since \eqref{critical} holds, the function $u_0\in H^1(\R^d)$ solves $(S_\lambda)$. \hfill$\Box$

\subsection{Some test functions with symmetries}\label{testf} Following Krist\'{a}ly, Moro\c{s}anu, and O'Regan \cite{KRO}, we construct some special test functions belonging to $Fix_{O(d)}(H^1(\R^d))$ that will be useful for our purposes. If $a<b$, define
$$
A_a^b:=\{x\in \R^d:a\leq |x|\leq b\}.
$$
Since $W\in L^{\infty}(\R^d)\setminus\{0\}$ is a radially symmetric function with $W\geq 0$, one can find real numbers $R>r>0$ and $\alpha>0$ such that
\begin{equation}\label{ess}
\displaystyle{\rm essinf}_{x\in A_{r}^{R}}W(x)\geq \alpha>0.
\end{equation}
\indent Hence, let $0<r<R$, such that \eqref{ess} holds and take $\sigma\in (0,(R-r)/2)$. Set $v_\sigma\in Fix_{O(d)}(H^1(\R^d))$ given by
\[
v_\sigma(x):= \left\{
\begin{array}{ll}
\displaystyle {\left(\frac{|x|-r}{\sigma}\right)_+} & \mbox{ if $|x|\leq r+\sigma$} \\
\displaystyle 1 & \mbox{ if $r+\sigma\leq |x|\leq R-\sigma$} \\
\displaystyle {\left(\frac{R-|x|}{\sigma}\right)_+} & \mbox{ if $|x|\geq R-\sigma$}
\end{array}
\right.
\]
where $z_+:=\max\{0,z\}$.
With the above notation, we have:
 \begin{itemize}
 \item[$(i_1)$] $\textrm{supp}(v_\sigma)\subseteq A_r^R$;
 \item[$(i_2)$] $\|v_\sigma\|_\infty\leq 1$;
 \item[$(i_3)$] $v_\sigma(x)=1$ for every $x\in A_{r+\sigma}^{R-\sigma}$.
 \end{itemize}

Now, assume $r\geq \displaystyle\frac{R}{5+4\sqrt{2}}$ and set $\sigma\in (0,1)$. Define $v_\sigma^{i}\in H^1(\R^d)$ as follows\par

\[
v^i_\sigma(x):= \left\{
\begin{array}{ll}
\displaystyle v^{\frac{d-2}{2}}_\sigma(x) & \mbox{ if $i=\displaystyle\frac{d-2}{2}$ and $x:=(x_1,x_3)\in \R^{d/2}\times\R^{d/2}$} \\
\displaystyle v^{\sigma}_i(x) & \mbox{ if $i\neq\displaystyle\frac{d-2}{2}$ and $x:=(x_1,x_2,x_3)\in \R^{i+1}\times \R^{d-2i-2}\times\R^{i+1}$},
\end{array}
\right.
\]
\noindent for every $x\in \R^d$, where:\par
$
\displaystyle v_\sigma^{\frac{d-2}{2}}(x_1,x_3):= \Bigg[\Bigg(\frac{R-r}{4}-\max\left\{\sqrt{\left(|x_1|^2-\frac{R+3r}{4}\right)^{2}+|x_3|^2}, \sigma\frac{R-r}{4}\right\}\Bigg)_+
$
$$
\quad\quad\quad
-\Bigg(\frac{R-r}{4}-\max\left\{\sqrt{\left(|x_1|^2-\frac{R+3r}{4}\right)^{2}+|x_3|^2}, \sigma\frac{R-r}{4}\right\}\Bigg)_+\Bigg]
$$
$$
\quad\quad\quad\times\frac{4}{(R-r)(1-\sigma)},\quad\quad\forall\, (x_1,x_3)\in \R^{d/2}\times\R^{d/2},
$$
and\par
$
\displaystyle v_{i}^{\sigma}(x_1,x_2,x_3):= \Bigg[\Bigg(\frac{R-r}{4}-\max\left\{\sqrt{\left(|x_1|^2-\frac{R+3r}{4}\right)^{2}+|x_3|^2}, \sigma\frac{R-r}{4}\right\}\Bigg)_+
$
$$
\quad\quad\quad
-\Bigg(\frac{R-r}{4}-\max\left\{\sqrt{\left(|x_3|^2-\frac{R+3r}{4}\right)^{2}+|x_1|^2}, \sigma\frac{R-r}{4}\right\}\Bigg)_+\Bigg]
$$
$$
\times\left(\frac{R-r}{4}-\max\left\{|x_2|,\sigma\frac{R-r}{4}\right\}\right)_+\frac{4}{(R-r)^2(1-\sigma)^2},
$$
for every $(x_1,x_2,x_3)\in \R^{d/2}\times\R^{d-2i-2}\times\R^{d/2}$, and $i\neq\displaystyle\frac{d-2}{2}$.\par
 Now, it is possible to prove that $v_\sigma^i\in Fix_{H_{d,\eta_i}}(H^1(\R^d))$. Moreover, for every $\sigma\in (0,1]$, let
 $$
 Q^{(1)}_\sigma:=\left\{(x_1,x_3)\in \R^{i+1}\times\R^{i+1}: \sqrt{\left(|x_1|^2-\frac{R+3r}{4}\right)^{2}+|x_3|^2}\leq\sigma\frac{R-r}{4}\right\}
 $$
 and
 $$
 Q^{(2)}_\sigma:=\left\{(x_1,x_3)\in \R^{i+1}\times\R^{i+1}: \sqrt{\left(|x_3|^2-\frac{R+3r}{4}\right)^{2}+|x_1|^2}\leq\sigma\frac{R-r}{4}\right\}.
 $$
\indent Define
 \[
D^i_\sigma:= \left\{
\begin{array}{ll}
\displaystyle D^{\frac{d-2}{2}}_\sigma & \mbox{ if $i=\displaystyle\frac{d-2}{2}$}\\
\displaystyle D^{\sigma}_i & \mbox{ if $i\neq\displaystyle\frac{d-2}{2}$},
\end{array}
\right.
\]
 where
 $$
 D_\sigma^{\frac{d-2}{2}}:=\left\{(x_1,x_3)\in \R^{d/2}\times\R^{d/2}: (x_1,x_3)\in Q^{(1)}_\sigma\cap Q^{(2)}_\sigma\right\},
 $$
 and
 $$
 D^{\sigma}_i:=\left\{(x_1,x_2,x_3)\in \R^{d/2}\times\R^{d-2i-2}\times\R^{d/2}: (x_1,x_3)\in Q^{(1)}_\sigma\cap Q^{(2)}_\sigma,\,\,{\rm and}\,\,|x_2|\leq \sigma\frac{R-r}{4}\right\},
 $$
for every $i\neq \displaystyle\frac{d-2}{2}$.\par
 \indent The sets $D^{i}_\sigma$ have positive Lebesgue measure and they are $H_{d,\eta_i}$-invariant. Moreover, for every $\sigma\in (0,1)$, one has $v_\sigma^i\in Fix_{H_{d,\eta_i}}(H^1(\R^d))$ and the following facts hold:
 \begin{itemize}
 \item[$(j_1)$] $\textrm{supp}(v_\sigma^i)=D_1^i\subseteq A[r,R]$;
 \item[$(j_2)$] $\|v_\sigma^i\|_\infty\leq 1$;
 \item[$(j_3)$] $|v_\sigma^i(x)|=1$ for every $x\in D_\sigma^i$.
 \end{itemize}

\section{Proof of the Main Result}\label{sec:Main}

Part $(a_1)$ - The main idea of the proof consists of applying Theorem \ref{BMB} to
the functional
$$\mathcal J_{\lambda}(u)= \Phi(u)-\lambda \Psi|_{Fix_{O(d)}(H^1(\R^d))}(u),\,\,\,\,\forall\, u \in  Fix_{O(d)}(H^1(\R^d)),$$
with
$$\Phi(u):=\frac 1 2 \left(\int_{\R^d}|\nabla u(x)|^2dx+\int_{\R^d}|u(x)|^2dx\right),$$
as well as
$$\Psi(u):= \int_{\mathbb R^d}W(x)F(u(x))dx.$$

Successively, the existence of one
non-trivial radial solution of problem $(S_\lambda)$ follows by the symmetric criticality principle due to Krawcewicz and Marzantowicz and recalled above, in Theorem \ref{Palais}.

To this aim, first notice that the functionals~$\Phi$ and
$\Psi|_{Fix_{O(d)}(H^1(\R^d))}$ have the regularity required by Theorem \ref{BMB}, according to Corollary~\ref{semicontinuity}. On the other hand, the functional $\Phi$ is clearly coercive in $Fix_{O(d)}(H^1(\R^d))$ and
$$\displaystyle\inf_{u\in Fix_{O(d)}(H^1(\R^d))}\Phi(u)=0.$$

Now, let us define
\begin{equation}\label{lambda1}
\lambda^{\star}:=\frac{1}{\kappa_1c_q}\max_{\gamma>0}\left(\frac{\gamma}{\displaystyle \sqrt{2}\|W\|_{\frac{q}{q-1}}+2^{q/2}c_q^{q-1}\|W\|_{\infty}\gamma^{q-1}}\right),
\end{equation}
where $\kappa_1=$ and
$$
c_\ell:=\sup\left\{\frac{\|u\|_\ell}{\|u\|}:u\in Fix_{O(d)}(H^1(\R^d))\setminus\{0\}\right\},
$$
for every $q\in (2,2^*)$ and take $0<\lambda<\lambda^{\star}$.

 Thanks to \eqref{lambda1},
there exists $\bar{\gamma}>0$ such that
\begin{equation}\label{n3}
\lambda<\lambda^{\star}{(\bar{\gamma})}:=\frac{{\bar\gamma}}{\kappa_1c_q}\left(\frac{1}{\displaystyle \sqrt{2}\|W\|_{\frac{q}{q-1}}+2^{q/2}c_q^{q-1}\|W\|_{\infty}\bar\gamma^{q-1}}\right).
\end{equation}
\indent Arguing as in \cite{Molica}, let us define the function $\chi: (0,+\infty)\rightarrow [0,+\infty)$ as
$$
\chi(r):=\frac{\displaystyle\sup_{u\in\Phi^{-1}((-\infty,r))}\Psi|_{Fix_{O(d)}(H^1(\R^d))}(u)}{r},
$$
for every $r>0$.\par
\indent It follows by \eqref{ineq2}  that
\begin{equation}\label{nghl}
\Psi|_{Fix_{O(d)}(H^1(\R^d))}(u) \leq \kappa_1\left(\|W\|_{\frac{q}{q-1}}\|u\|_q+\|W\|_{\infty}\|u\|_q^q\right),
\end{equation}
for every $u\in Fix_{O(d)}(H^1(\R^d))$.\par
\indent Moreover, one has
\begin{equation}\label{gio}
\|u\|< \sqrt{2r},
\end{equation}
\noindent for every $u\in\Phi^{-1}((-\infty,r))$.\par
\indent Now, by using (\ref{gio}), the Sobolev embedding \eqref{compattezza} and \eqref{nghl} yield
$$
\Psi|_{Fix_{O(d)}(H^1(\R^d))}(u)< \kappa_1c_q\left(\|W\|_{\frac{q}{q-1}}\sqrt{{2r}}+c_q^{q-1}\|W\|_{\infty}(2r)^{q/2}\right),
$$
for every $u\in\Phi^{-1}((-\infty,r))$.\par
\noindent Consequently,
$$
\sup_{u\in\Phi^{-1}((-\infty,r))}\Psi|_{Fix_{O(d)}(H^1(\R^d))}(u)\leq \kappa_1 c_q\left(\|W\|_{\frac{q}{q-1}}\sqrt{{2r}}+c_q^{q-1}\|W\|_{\infty}(2r)^{q/2}\right).
$$
\indent The above inequality yields
\begin{equation}\label{n}
\chi(r)\leq \kappa_1 c_q\left(\|W\|_{\frac{q}{q-1}}\sqrt{\frac{2}{r}}+2^{q/2}c_q^{q-1}\|W\|_{\infty}r^{q/2-1}\right),
\end{equation}
\noindent for every $r>0$.\par
\indent Evaluating inequality \eqref{n} in $r=\bar\gamma^2$, it follows that
\begin{equation}\label{n2}
\chi(\bar{\gamma}^2)\leq \kappa_1 c_q\left(\sqrt{{2}}\frac{\|W\|_{\frac{q}{q-1}}}{\bar\gamma}+2^{q/2}c_q^{q-1}\|W\|_{\infty}{\bar\gamma}^{q-2}\right).
\end{equation}
Now, we notice that
$$
\varphi(\bar{\gamma}^2):=\inf_{u\in\Phi^{-1}((-\infty,\bar{\gamma}^2))}\frac{\displaystyle\left(\sup_{v\in\Phi^{-1}((-\infty,\bar{\gamma}^2))}\Psi|_{Fix_{O(d)}(H^1(\R^d))}(v)\right)-\Psi|_{Fix_{O(d)}(H^1(\R^d))}(u)}{r-\Phi(u)}\leq \chi(\bar{\gamma}^2),
$$
owing to $z_0\in\Phi^{-1}((-\infty,\bar{\gamma}^2))$ and $\Phi(z_0)=\Psi|_{Fix_{O(d)}(H^1(\R^d))}(z_0)=0$, where $z_0\in Fix_{O(d)}(H^1(\R^d))$ is the zero function.\par
 Thanks to (\ref{n3}), the above inequality in addition to (\ref{n2}) give
\begin{equation}\label{n2fg}
\varphi(\bar{\gamma}^2)\leq \chi(\bar{\gamma}^2)\leq \kappa_1c_q\left(\sqrt{{2}}\frac{\|W\|_{\frac{q}{q-1}}}{\bar\gamma}+2^{q/2}c_q^{q-1}\|W\|_{\infty}{\bar\gamma}^{q-2}\right)<\frac{1}{\lambda}.
\end{equation}
In conclusion,
$$
\lambda\in \left(0,\frac{{\bar\gamma}}{\kappa_1c_q}\left(\frac{1}{\displaystyle \sqrt{2}\|W\|_{\frac{q}{q-1}}+2^{q/2}c_q^{q-1}\|W\|_{\infty}\bar\gamma^{q-1}}\right)\right)\subseteq (0,{1}/{\varphi(\bar{\gamma}^2)}).
$$

Invoking Theorem \ref{BMB}, there exists a function $u_\lambda\in\Phi^{-1}((-\infty,\bar{\gamma}^2))$ such that
$$\mathcal{J}^0(u_\lambda;\varphi)\geq 0,\quad \forall\,\varphi\in Fix_{O(d)}(H^1(\R^d)).$$
More precisely, the function $u_\lambda$ is a global minimum of the restriction of the functional $\mathcal{J}_{\lambda}$ to the sublevel $\Phi^{-1}((-\infty,\bar{\gamma}^2))$.

 Hence,
let $u_\lambda$ be such that
\begin{equation}\label{minimum}
\mathcal J_{\lambda}(u_{\lambda})\leq \mathcal
J_{\lambda}(u),\quad \mbox{for any}\,\,\, u\in Fix_{O(d)}(H^1(\R^d))\,\,
\mbox{such that}\,\,\,\Phi(u)<\bar{\gamma}^2
\end{equation}
and
\begin{equation}\label{minimum2}
\Phi(u_{\lambda})<\bar\gamma^2\,,
\end{equation}
and also $u_\lambda$ is  a critical point of $\mathcal
J_{\lambda}$ in $Fix_{O(d)}(H^1(\R^d))$.
Now, the orthogonal group $O(d)$ acts isometrically on $H^1(\R^d)$ and, thanks to the symmetry of the potential $W$, one has
$$
\int_{\R^d}W(x)F((gu)(x))dx=\int_{\R^d}W(x)F(u(g^{-1}x))dx=\int_{\R^d}W(z)F(u(z))dz,
$$
for every $g\in O(d)$. Then the functional $J_\lambda$ is $O(d)$-invariant on $H^1(\R^d)$.\par
  So, owing to Theorem \ref{Palais}, $u_\lambda$ is a weak solution of
problem~$(S_\lambda)$\,.
In this
setting, in order to prove that $u_\lambda\not \equiv 0$ in $Fix_{O(d)}(H^1(\R^d))$\,,
first we claim that there exists a sequence of functions $\big\{w_j\big\}_{j\in
\NN}$ in $Fix_{O(d)}(H^1(\R^d))$ such that
\begin{equation}\label{wjBlu}
\limsup_{j\to
+\infty}\frac{\Psi|_{Fix_{O(d)}(H^1(\R^d))}(w_j)}{\Phi(w_j)}=+\infty\,.
\end{equation}

By the assumption on the limsup in~(\ref{azero}), there exists a sequence
$\{s_j\}_{j\in\NN}\subset (0,+\infty)$ such that $s_j\to 0^+$ as $j\to
+\infty$ and
\begin{equation}\label{limsup}
 \lim_{j\rightarrow +\infty}\frac{\displaystyle F(s_j)}{s_j^2}=+\infty,
\end{equation}
namely, we have that for any $M>0$ and $j$
sufficiently large
\begin{equation}\label{M}
\displaystyle\,F(s_j)>Ms_j^2\,.
\end{equation}

\indent Now, define $w_j:=s_jv_\sigma$ for any $j\in \NN$, where the function $v_\sigma$ is given in Subsection \ref{testf}. Since $v_\sigma\in Fix_{O(d)}(H^1(\R^d))$ of course, one has $w_j\in Fix_{O(d)}(H^1(\R^d))$ for any $j\in \NN$. Bearing in mind that the functions $v_\sigma$ satisfy  $(i_1)$--$(i_3)$, thanks to $F(0)=0$ and (\ref{M}) we have
\begin{eqnarray}\label{conto1}
& {\frac{\displaystyle \Psi|_{Fix_{O(d)}(H^1(\R^d))}(w_j)}{\displaystyle\Phi(w_j )}}= \frac{\displaystyle \int_{A_{r+\sigma}^{R-\sigma}}W(x)F(w_j(x))\, dx+\int_{A_{r}^{R}\setminus
A_{r+\sigma}^{R-\sigma}} W(x)F(w_j(x))\,dx}{\displaystyle \Phi(w_j)} \nonumber \\
& \qquad \qquad= \frac{\displaystyle \int_{A_{r+\sigma}^{R-\sigma}}W(x)F(s_j)\, dx+\int_{A_{r}^{R}\setminus A_{r+\sigma}^{R-\sigma}}
W(x)F(s_jv_\sigma(x))\,dx}{\displaystyle\Phi(w_j)}  \\
& \qquad \qquad\quad\geq \displaystyle 2\frac{M|A_{r+\sigma}^{R-\sigma}|\alpha s^2_j+{\displaystyle \int_{A_{r}^{R}\setminus A_{r+\sigma}^{R-\sigma}}
W(x)F(s_jv_\sigma(x))\,dx}}{s_j^2\|v_\sigma\|^2}\,,\nonumber
\end{eqnarray}
for $j$ sufficiently large.

Now, we have to consider two different cases.\par

{Case 1:} ${\displaystyle \lim_{s\to 0^+}\frac{\displaystyle F(s)}{s^2}=+\infty}$.

Then there exists $\rho_M>0$ such that for any $s$ with $0<s<\rho_M$
\begin{equation}\label{H}
F(s) \geq Ms^2\,.
\end{equation}

Since $s_j\to 0^+$ and $0\leq v_\sigma(x)\leq 1$ in $\R^d$, it follows that $w_j(x)=s_jv_\sigma(x)\to 0^+$ as $j\to +\infty$ uniformly in $x\in \R^d$. Hence, $0\leq w_j(x)< \rho_M$ for $j$ sufficiently large and for any $x\in \R^d$. Hence, as a consequence of (\ref{conto1}) and (\ref{H}), we have
that
\begin{eqnarray}
& {\displaystyle \frac{\Psi|_{Fix_{O(d)}(H^1(\R^d))}(w_j)}{\Phi(w_j )}  \geq 2\frac{M|A_{r+\sigma}^{R-\sigma}|\alpha
s_j^2+{\displaystyle \int_{A_{r}^{R}\setminus
A_{r+\sigma}^{R-\sigma}} W(x)F(s_jv_\sigma(x))\,dx}}{s_j^2\|v_\sigma\|^2}} \nonumber\\
& {\displaystyle \qquad \qquad  \geq 2M\alpha\frac{|A_{r+\sigma}^{R-\sigma}|+{\displaystyle \int_{A_{r}^{R}\setminus
A_{r+\sigma}^{R-\sigma}}
|v_\sigma(x)|^2\,dx}}{\|v_\sigma\|^2}}\,, \nonumber
\end{eqnarray}
for $j$ sufficiently large. The arbitrariness of $M$ gives
(\ref{wjBlu}) and so the claim is proved.

\medskip

{Case 2:} ${\displaystyle \liminf_{s\to 0^+}\frac{\displaystyle F(s)}{s^2}=\ell\in
\RR}$\,.

Then for any $\varepsilon>0$ there exists $\rho_\varepsilon>0$ such that for any $s$ with $0<s<\rho_\varepsilon$
\begin{equation}\label{epsilon}
F(s) \geq(\ell-\varepsilon)s^2\,.
\end{equation}
Arguing as above, we can suppose that $0\leq w_j(x)=s_jv_\sigma(x)<\rho_\varepsilon$ for $j$ large enough and any $x\in \R^d$. Thus, by (\ref{conto1}) and (\ref{epsilon}) we
get
\begin{eqnarray}\label{eps2}
& {\displaystyle \frac{\Psi|_{Fix_{O(d)}(H^1(\R^d))}(w_j)}{\Phi(w_j )}  \geq
2\frac{M |A_{r+\sigma}^{R-\sigma}|\alpha s_j^2+{\displaystyle \int_{A_{r}^{R}\setminus
A_{r+\sigma}^{R-\sigma}} W(x)F(
s_jv_\sigma(x))\,dx}}{s_j^2\|v_\sigma\|^2}}\\
&\qquad \qquad \qquad\quad {\displaystyle \geq 2\alpha\frac{
M|A_{r+\sigma}^{R-\sigma}|+{(\ell-\varepsilon)\displaystyle
\int_{A_{r}^{R}\setminus
A_{r+\sigma}^{R-\sigma}} |v_\sigma(x)|^2\,dx}}{\|v_\sigma\|^2}},\nonumber
\end{eqnarray}
provided that $j$ is sufficiently large.

Let
$$M>\max\left\{0, -\frac{2\ell}{|A_{r+\sigma}^{R-\sigma}|}\int_{A_{r}^{R}\setminus
A_{r+\sigma}^{R-\sigma}}|v_\sigma(x)|^2\,dx\right\}\,,$$
and
$$0<\varepsilon <\frac{\displaystyle\frac{M}{2}|A_{r+\sigma}^{R-\sigma}|+\ell \int_{A_{r}^{R}\setminus
A_{r+\sigma}^{R-\sigma}}|v_\sigma(x)|^2\,dx}{\displaystyle\int_{A_{r}^{R}\setminus
A_{r+\sigma}^{R-\sigma}}|v_\sigma(x)|^2\,dx}.$$

By (\ref{eps2}) we have
\begin{eqnarray}
& {\displaystyle \frac{\Psi|_{Fix_{O(d)}(H^1(\R^d))}(w_j)}{\Phi(w_j )}  \geq
2\alpha\frac{
M|A_{r+\sigma}^{R-\sigma}|+{(\ell-\varepsilon)\displaystyle
\int_{A_{r}^{R}\setminus
A_{r+\sigma}^{R-\sigma}} |v_\sigma(x)|^2\,dx}}{\|v_\sigma\|^2}} \nonumber \\
& \qquad\qquad\quad \displaystyle \geq
\frac{2\alpha}{\|v_\sigma\|^2}\left(M|A_{r+\sigma}^{R-\sigma}|+{\ell\displaystyle
\int_{A_{r}^{R}\setminus
A_{r+\sigma}^{R-\sigma}} |v_\sigma(x)|^2\,dx}-\varepsilon {\displaystyle\int_{A_{r}^{R}\setminus
A_{r+\sigma}^{R-\sigma}}|v_\sigma(x)|^2\,dx}\right)\nonumber\\
&  {\displaystyle \geq\alpha M \frac{|A_{r+\sigma}^{R-\sigma}|}{\|v_\sigma\|^2}},\nonumber
 \end{eqnarray}
for $j$ sufficiently large. Hence, assertion~(\ref{wjBlu}) is clearly verified.

Now, we notice that
$$\|w_j\|=s_j\,\|v_\sigma\|\to 0,$$
as $j\to +\infty$\,, so that for $j$ large enough
$$\|w_j\|< \sqrt{2}\bar \gamma.$$
\indent Hence
\begin{equation}\label{wj1}
w_j\in \Phi^{-1}\big((-\infty, \bar \gamma^2)\big)\,,
\end{equation}
and on account of \eqref{wjBlu}, also
\begin{equation}\label{wj2}
\mathcal
J_{\lambda}(w_j)=\Phi(w_j)-\lambda\Psi|_{Fix_{O(d)}(H^1(\R^d))}(w_j)<0,
\end{equation}
for $j$ sufficiently large.

Since $u_\lambda$ is a global minimum of the restriction $\mathcal
J_{\lambda}|_{\Phi^{-1}((-\infty,\bar \gamma^2))}$, by (\ref{wj1}) and (\ref{wj2})
we have that
\begin{equation}\label{nontrivial}
\mathcal J_{\lambda}(u_\lambda)\leq \mathcal
J_{\lambda}(w_j)<0=\mathcal J_{\lambda}(0)\,,
\end{equation}
so that $u_\lambda\not\equiv 0$ in $Fix_{O(d)}(H^1(\R^d))$.\par
 Thus, $u_\lambda$ is a
non-trivial weak solution of problem~$(S_\lambda)$. The arbitrariness
of $\lambda$ gives that $u_\lambda\not \equiv 0$ for any
$\lambda\in (0, \lambda^{\star})$. By a Strauss-type estimate (see Lions \cite{Lions}) we have that
$|u_\lambda(x)|\rightarrow 0$ as $|x|\rightarrow \infty$.
This concludes the proof of part $(a_1)$ of Theorem~\ref{Main1}.\par

\smallskip

Part $(a_2)$ - Let
$$
c_{i,\ell}:=\sup\left\{\frac{\|u\|_\ell}{\|u\|}:u\in Fix_{H_{d,\eta_i}}(H^1(\R^d))\setminus\{0\}\right\},
$$
for every $\ell\in (2,2^*)$, with $i\in J_d$ and set
\begin{equation}\label{lambda2}
\lambda^{\star}_{i,q}:=\frac{1}{\kappa_1 c_{i,q}}\max_{\gamma>0}\left(\frac{{\gamma}}{\displaystyle \sqrt{2}\|W\|_{\frac{q}{q-1}}+2^{q/2}c_{i,q}^{q-1}\|W\|_{\infty}\gamma^{q-1}}\right).
\end{equation}
Assume $d> 3$ and suppose that the potential $F$ is even. Let
$$
\lambda_{\star}:= \left\{
\begin{array}{l}
\lambda^{\star} \,\,\,\,{\rm if}\,\,\,d=5\\
\smallskip
\min\{\lambda^{\star},\lambda^{\star}_{i,q}:i\in J_d\}\,\,\,\,{\rm if}\,\,\,d\neq 5.
\end{array}
\right.
$$
We claim that for every $\lambda\in (0,\lambda_{\star})$ problem $(S_\lambda)$
 admits at least
 $$
 \zeta^{(d)}_S:=1+(-1)^{d}+\left[\frac{d-3}{2}\right]
 $$
 pairs of non-trivial weak solutions $\{\pm u_{\lambda,i}\}_{i\in J'_d}\subset H^1(\R^d)$, where $J'_d:=\{1,...,\zeta^{(d)}_S\}$, such that
$|u_{\lambda,i}(x)|\rightarrow 0$, as $|x|\rightarrow \infty$, for every $i\in J'_d$. \par
  Moreover, if $d\neq 5$ problem $(S_\lambda)$ admits at least
 $$\tau_d:=(-1)^{d}+\left[\frac{d-3}{2}\right]$$ pairs of sign-changing weak solutions.\par
We divide the proof into two parts.\par
{Part 1:} dimension $d=5$. Since $F$ is symmetric, the energy functional
$$\mathcal J_{\lambda}(u):= \Phi(u)-\lambda \Psi|_{Fix_{O(d)}(H^1(\R^d))}(u),\,\,\,\,\forall\, u \in  Fix_{O(d)}(H^1(\R^d)),$$
is even. Owing to Theorem \ref{Main1}, for every $\lambda\in (0,\lambda^{\star})$, problem $(S_\lambda)$ admits at least one (that is $\zeta^{(5)}_S=1$) non-trivial pair of radial weak solutions $\{\pm u_{\lambda}\}\subset H^1(\R^d)$. Furthermore, the functions $\pm u_\lambda$ are homoclinic.\par
\smallskip
{Part 2:} dimension $d>3$ and $d\neq 5$.
For every $\lambda>0$ and $i\in J_d$, consider the restriction
$\mathcal H_{\lambda, i}:=J_\lambda|_{Fix_{H_{d,\eta_i}}(H^1(\R^d))}:Fix_{H_{d,\eta_i}}(H^1(\R^d))\to \mathbb R$ defined by
$$\mathcal H_{\lambda, i}:=\Phi_{H_{d,\eta_i}}(u)-\lambda \Psi|_{Fix_{H_{d,\eta_i}}(H^1(\R^d))}(u),$$ where
\begin{equation*}
    \Phi_{H_{d,\eta_i}}(u):=\frac{1}{2}\|u\|^2\,\,\, {\rm and}\,\,\, \Psi|_{Fix_{H_{d,\eta_i}}(H^1(\R^d))}(u):=\int_{\mathbb R^d}W(x)F(u(x))dx,
\end{equation*}
for every $u\in Fix_{H_{d,\eta_i}}(H^1(\R^d))$.\par
In order to obtain the existence of
$$\tau_d:=(-1)^{d}+\left[\frac{d-3}{2}\right]$$ pairs of sign-changing weak solutions $\{\pm z_{\lambda,i}\}_{i\in J_d}\subset H^1(\R^d)$, where $J_d:=\{1,...,\tau_d\}$,
 the main idea of the proof consists in applying Theorem \ref{BMB} to
the functionals~$\mathcal H_{\lambda,i}$, for every $i\in J_d$. We notice that, since $d>3$ and $d\neq 5$, $\tau_d\geq 1$. Consequently, the cardinality $|J_d|\geq 1$.

\indent Since $0<\lambda<\lambda_{i,q}^{\star}$, with $i\in J_d$,
there exists $\bar{\gamma}_i>0$ such that
\begin{equation}\label{n32a}
\lambda<\lambda_{\star}^{(i)}{(\bar{\gamma}_i)}:=\frac{{\bar\gamma_i}}{\kappa_1 c_{i,q}}\left(\frac{1}{\displaystyle \sqrt{2}\|W\|_{\frac{q}{q-1}}+2^{q/2}c_{i,q}^{q-1}\|W\|_{\infty}\bar\gamma^{q-1}_i}\right).
\end{equation}
Similar arguments used for proving \eqref{n2fg} yield
\begin{equation}\label{n2fgnew}
\varphi(\bar{\gamma}^2_i)\leq \chi(\bar{\gamma}^2_i)\leq \kappa_1c_q\left(\sqrt{{2}}\frac{\|W\|_{\frac{q}{q-1}}}{\bar\gamma_i}+{2^{q/2}c_q^{q-1}}\|W\|_{\infty}{\bar\gamma_i}^{q-2}\right)<\frac{1}{\lambda}.
\end{equation}
\indent Thus,
$$
\lambda\in \left(0,\frac{{\bar\gamma}_i}{\kappa_1c_q}\left(\frac{1}{\displaystyle \sqrt{2}\|W\|_{\frac{q}{q-1}}+2^{q/2}c_q^{q-1}\|W\|_{\infty}\bar\gamma^{q-1}_i}\right)\right)\subseteq (0,{1}/{\varphi(\bar{\gamma}^2_i)}).
$$
Thanks to Theorem \ref{BMB}, there exists a function $z_{\lambda,i}\in \Phi_{H_{d,\eta_i}}^{-1}((-\infty,\bar{\gamma}^2_i))$ such that
$$\mathcal{J}^0(z_{\lambda,i};\varphi)\geq 0,\quad \forall\,\varphi\in Fix_{H_{d,\eta_i}}(H^1(\R^d))$$
\noindent and, in particular, $z_{\lambda,i}$ is a global minimum of the restriction of $\mathcal{H}_{\lambda, i}$ to $\Phi_{H_{d,\eta_i}}^{-1}((-\infty,\bar{\gamma}^2_i))$.

Due to the evenness of $J_\lambda$, bearing in mind \eqref{isometr}, and thanks to the symmetry assumptions on the potential $W$, we have
that the functional $J_{\lambda}$ is $H_{d,\eta_i}$-invariant on $H^1(\R^d)$, i.e.
$$J_{\lambda}(h\circledast_i u)=J_{\lambda}(u),$$  for every $h\in H_{d,\eta_i}$ and $u\in H^1(\R^d)$. Indeed, the group $H_{d,\eta_i}$ acts isometrically on $H^1(\R^d)$ and, thanks to the symmetry assumption on $W$, it follows that
$$
\int_{\R^d}W(x)F((hu)(x))dx=\int_{\R^d}W(x)F(u(h^{-1}x))dx=\int_{\R^d}W(z)F(u(z))dz,
$$
if $h\in H_{d,i}$, and
$$
\int_{\R^d}W(x)F((hu)(x))dx=\int_{\R^d}W(x)F(u(g^{-1}\eta_{_{H_{d,i}}}^{-1}x))dx=\int_{\R^d}W(z)F(u(z))dz,
$$
if $h=\eta_{_{H_{d,i}}} g\in H_{d,\eta_i}\setminus H_{d,i}$.\par
On account of Theorem \ref{Palais}, the critical point pairs $\{\pm z_{\lambda,i}\}$ of $\mathcal H_{\lambda,i}$ are also (generalized) critical points of $J_\lambda$.\par

Let $z_{\lambda,i}\in Fix_{H_{d,\eta_i}}(H^1(\R^d))$ be a critical point of $\mathcal
H_{\lambda,i}$ in $Fix_{H_{d,\eta_i}}(H^1(\R^d))$ such that
\begin{equation}\label{minimum23444}
\mathcal H_{\lambda,i}(z_{\lambda,i})\leq \mathcal
H_{\lambda,i}(u),\quad \mbox{for any}\,\,\, u\in Fix_{H_{d,\eta_i}}(H^1(\R^d))\,\,
\mbox{such that}\,\,\,\Phi_{H_{d,\eta_i}}(u)<\bar{\gamma}^2_i
\end{equation}
and
\begin{equation}\label{minimum24444}
\Phi_{H_{d,\eta_i}}(z_{\lambda,i})<\bar{\gamma}^2_i.
\end{equation}

In order to prove that $z_{\lambda,i}\not \equiv 0$ in $Fix_{H_{d,\eta_i}}(H^1(\R^d))$\,,
we claim that there exists a sequence $\big\{w_j^{i}\big\}_{j\in
\NN}$ in $Fix_{H_{d,\eta_i}}(H^1(\R^d))$ such that
\begin{equation}\label{wj2gg}
\limsup_{j\to
+\infty}\frac{\Psi|_{Fix_{H_{d,\eta_i}}(H^1(\R^d))}(w_j^i)}{\Phi(w_j^i)}=+\infty\,.
\end{equation}

\indent The sequence $\big\{w_j^{i}\big\}_{j\in
\NN} \subset Fix_{H_{d,\eta_i}}(H^1(\R^d))$, for which \eqref{wj2gg} holds, can be constructed by using the test functions introduced in \cite{KRO} and recalled in Subsection \ref{testf}. Thus, let us define $w_j^{i}:=s_jv_\sigma^i$ for any $j\in \NN$. Clearly, $w_j^{i}\in Fix_{H_{d,\eta_i}}(H^1(\R^d))$ for any $j\in \NN$. Moreover, taking into account the properties of $v_\sigma^i$ displayed in $(j_1)$--$(j_3)$, simple computations show that
\begin{eqnarray}\label{conto1344}
& {\frac{\displaystyle \Psi|_{Fix_{H_{d,\eta_i}}(H^1(\R^d))}(w_j^{i})}{\displaystyle\Phi(w_j^{i} )}}= \frac{\displaystyle \int_{D_\sigma^i}W(x)F(w_j^{i}(x))\, dx+\int_{A_{r}^{R}\setminus
D_\sigma^i} W(x)F(w_j^{i}(x))\,dx}{\displaystyle \Phi(w_j^{i})} \nonumber \\
& \qquad \qquad= \frac{\displaystyle \int_{D_\sigma^i}W(x)F(s_j)\, dx+\int_{A_{r}^{R}\setminus D_\sigma^i}
W(x)F(s_jv_\sigma^i(x))\,dx}{\displaystyle\Phi(w_j^{i})}  \\
& \qquad \qquad\quad\geq \displaystyle 2\frac{M|D_\sigma^i|\alpha s^2_j+{\displaystyle \int_{A_{r}^{R}\setminus D_\sigma^i}
W(x)F(s_jv_\sigma^i(x))\,dx}}{s_j^2\|v_{\sigma}^i\|^2}\,,\nonumber
\end{eqnarray}
for $j$ sufficiently large.\par
Arguing as in the proof of Theorem \ref{Main1}, inequality \eqref{conto1344} yields \eqref{wj2gg} and consequently,
we conclude that
\begin{equation*}\label{nontrivialnew}
\mathcal H_{\lambda,i}(z_{\lambda,i})\leq \mathcal
H_{\lambda,i}(w_j^{i})<0=\mathcal H_{\lambda,i}(0)\,,
\end{equation*}
so that $z_{\lambda,i}\not\equiv 0$ in $Fix_{H_{d,\eta_i}}(H^1(\R^d))$. In addition, $|z_{\lambda,i}(x)|\rightarrow 0$ as $|x|\rightarrow \infty$.\par
\indent On the other hand, since $\lambda<\lambda^{\star}$ and $F$ is even, Theorem \ref{Main1} and the principle of symmetric criticality (recalled in Theorem \ref{Palais}) ensure that problem $(S_\lambda)$ admits
at least one non-trivial pair of radial weak solutions $\{\pm u_\lambda\}\subset H^1(\R^d)$. Moreover, $|u_\lambda(x)|\rightarrow 0$ as $|x|\rightarrow \infty$.\par
In conclusion, since $\lambda<\lambda_\star$, there exist $\tau_d+1$ positive numbers $\bar\gamma$, $\bar\gamma_1$,...,$\bar\gamma_{\tau_d}$ such that
$$
\pm u_\lambda \in \Phi^{-1}((-\infty,\bar{\gamma}^2))\setminus\{0\}\subset Fix_{O(d)}(H^1(\R^d)),
$$
and
$$
\pm z_{\lambda,i}\in \Phi^{-1}_{H_{d,\eta_i}}((-\infty,\bar{\gamma}^2_i))\setminus\{0\}\subset Fix_{H_{d,\eta_i}}(H^1(\R^d)).
$$
Bearing in mind relations \eqref{intersezione} and \eqref{intersezione2} of Proposition \ref{geometry} (see also \cite[Theorem 2.2]{KRO} for details) we have that $$\Phi^{-1}((-\infty,\bar{\gamma}^2))\cap \Phi^{-1}_{H_{d,\eta_i}}((-\infty,\bar{\gamma}^2_i))\setminus\{0\}=\emptyset,$$
for every $i\in J_d$ and
$$\Phi^{-1}_{H_{d,\eta_i}}((-\infty,\bar{\gamma}^2_i))\cap \Phi^{-1}_{H_{d,\eta_j}}((-\infty,\bar{\gamma}^2_j))\setminus\{0\}=\emptyset,$$
for every $i,j\in J_d$ and $i\neq j$. Consequently problem $(S_\lambda)$
admits at least
 $$
 \zeta^{(d)}_S:=\tau_d+1,
 $$
 pairs of non-trivial weak solutions $\{\pm u_{\lambda,i}\}_{i\in J'_d}\subset H^1(\R^d)$, where $J'_d:=\{1,...,\zeta^{(d)}_S\}$, such that
$|u_{\lambda,i}(x)|\rightarrow 0$, as $|x|\rightarrow \infty$, for every $i\in J'_d$. Moreover, by construction, it follows that
 $$\tau_d:=(-1)^{d}+\left[\frac{d-3}{2}\right]$$ pairs of the attained solutions are sign-changing.\par
 The proof is now complete.\hfill$\Box$\par

 \section{Some applications}\label{sec finale}

A simple prototype of a function $F$ fulfilling the structural assumption \eqref{crescita} can be easily constructed as follows.
 Let $f:\R\rightarrow \R$ be a measurable function such that
 \begin{equation}\label{crescitafunz}
\sup_{s\in \R}\frac{|f(s)|}{1+|s|^{q-1}}<+\infty,
\end{equation}
for some $q\in \left(2,2^*\right)$. Furthermore, let $F$ be the potential defined by
$$
F(s):=\int_0^sf(t)dt,
$$
for every $s\in \R$. Of course $F$ is a Carath\'{e}odory function that is locally Lipschitz with $F(0)=0$. Since the growth condition \eqref{crescitafunz} is satisfied, $f$ is locally essentially bounded,
that is $f\in L^{\infty}_{\rm loc}(\R^d)$. Thus, invoking \cite[Proposition 1.7]{Mot1} it follows that
\begin{equation}\label{crescitafunz2}
\partial F(s)=[\underline{f}(s),\overline{f}(s)]
\end{equation}
where
$$
\underline{f}(s):=\lim_{\delta\rightarrow 0^+}\displaystyle{\rm essinf}_{|t-s|<\delta}{f(t)},
$$
and
$$
\overline{f}(s):=\lim_{\delta\rightarrow 0^+}\displaystyle{\rm esssup}_{|t-s|<\delta}{f(t)},
$$
for every $s\in \R$.\par
On account of \eqref{crescitafunz} and \eqref{crescitafunz2}, inequality \eqref{crescita} immediately follows. Furthermore, if $f$ is a continuous function and \eqref{crescitafunz} holds, then problem $(S_\lambda)$ assumes the simple and significative form:
\begin{itemize}
\item[$({S}'_\lambda)$] {\it Find $u\in H^1(\R^d)$ such that
\begin{equation*}
\left\{\begin{array}{l} \displaystyle\int_{\R^d}\nabla u(x)\cdot \nabla \varphi(x)dx +\int_{\R^d}u(x)\varphi(x))dx\\
\qquad \qquad\qquad \qquad\quad\quad -\displaystyle\lambda\int_{\R^d}W(x)f(u(x))\varphi(x)dx= 0,\,\,\,\,\,\,\\
\forall\,\varphi\in H^1(\R^d).
\end{array}\right.
\end{equation*}}
\end{itemize}
See \cite{GR} for related topics.

Of course, the solutions of $({S}'_\lambda)$ are exactly the weak solutions of the following Schr\"{o}dinger equation
\begin{equation*}\label{Dirichlet}
\left\{
\begin{array}{l}
-\Delta u+u= \lambda W(x)f(u) \,\,\,\,{\rm in}\,\,\,\R^d\\
\smallskip
u\in H^{1}(\R^d),
\end{array}
\right.
\end{equation*}
which has been
widely studied in the literature. In particular, Theorem \ref{Main1} can be viewed as a non-smooth version of the results contained in \cite{Molica}. See, among others, the papers \cite{AM2, AM, Bartolo, BCW(MatANN), BW} as well as \cite{ClappWeth, Filippucci, MMP, seminalrabinowitz}.\par

We point out that the approach adopted here can be used in order to study the existence of multiple solutions for
hemivariational inequalities on a strip-like domain of the Euclidean space (see \cite{Kri2} for related topics).
Since this approach differs to the above, we will treat it in a forthcoming
paper.

\smallskip
\noindent {\bf Acknowledgements.} This research was realized under the auspices of
the Italian MIUR project \textit{Variational methods, with applications to
problems in mathematical physics and geometry} (2015KB9WPT 009) and the Slovenian Research Agency grants P1-0292, J1-8131, J1-7025, N1-0083, and N1-0064. \par

\end{document}